\documentclass[11pt,a4paper]{amsart}
\usepackage{amsmath,amsfonts,amsthm,amssymb,amscd}
\usepackage[latin1]{inputenc}
\usepackage{mathrsfs}
\usepackage{latexsym}
\usepackage{graphicx}
\usepackage{enumerate}
\def\M{\mathbb{N}_0}
\def\N{\mathbb{N}}
\def\W{{\rm W}}

\def\B{\mathcal{B}}
\def\P{\mathcal{P}}

\def\p2{\mathcal{P}(\W^2)}

\def\r{\mathcal{R}}
\def\T{\mathcal{T}}

\newtheorem{Thm}{Theorem}[section]
\newtheorem{Prop}[Thm]{Proposition}

\newtheorem{Def}{Definition}

\newtheorem{Rem}{Remark}
\newtheorem{Not}{Notation}
\begin{document}
\title{Blocking Wythoff Nim}

\author{Urban Larsson}
\email{urban.larsson@chalmers.se} 
\address{Mathematical Sciences, Chalmers University of Technology 
and University of Gothenburg, G\"oteborg, Sweden}
\keywords{Beatty sequence, Blocking maneuver, 
exact $k$-covers, Impartial game, Muller Twist, 
Wythoff Nim.}
\date{\today }

\begin{abstract}
The 2-player impartial game of Wythoff Nim is played on two piles of tokens. A move consists in removing any number of tokens from precisely one of the piles or the same number of tokens from both piles. The winner is the player who removes the last token. We study this game with a blocking maneuver, that is, for each move, before the next player moves the previous player may declare at most a predetermined number, $k - 1 \ge 0$, of the options as forbidden. When the next player has moved, any blocking maneuver is forgotten and does not have any further impact on the game. We resolve the winning strategy of this game for $k = 2$ and $k = 3$ and, supported by computer simulations, state conjectures of the asymptotic `behavior' of the $P$-positions for the respective games when $4 \le k \le 20$.
\end{abstract}

\maketitle
\vskip 30pt

\section{Introduction}\label{Sec:1}
In this note we study a variation of the 2-player 
combinatorial game of Wythoff Nim \cite{Wyt07}. 
The game is \emph{impartial}, because given a position in the game, the 
set of options does not depend on which player is in turn to move. A 
background on impartial games may be found in \cite{AlNoWo07, BeCoGu82, Con76}.
Let $\N$ and $\M$ denote the positive and non-negative integers 
respectively. Let the \emph{game board} be $\B:=\M\times \M$.
 
\begin{Def}\label{Def:1}
Let $(x, y)\in \B$. Then $(x - i, y - j)$ is an option  
of Wythoff Nim if either: 
\begin{itemize}
\item [(v)] $0 = i < j\le y$,
\item [(h)] $0 = j < i\le x$,
\item [(d)] $0 < i = j\le \min \{x, y\}$,
\end{itemize}
$i,j\in \N$.
\end{Def}
In this definition one might want to think about (v), (h) and (d) as 
symbolizing the `vertical' $(0, i)$ , `horizontal' $(i, 0)$ 
and `diagonal' $(i, i)$ \emph{moves} respectively. Two players 
take turns in moving according to these rules. The player who moves to the 
position $(0, 0)$ is declared the winner.
%Clearly, by these rules, the only terminal position of Wythoff Nim is $(0, 0)$. 
%The player who moves there is declared the winner. 
Here we study a variation of Wythoff Nim with 
a \emph{blocking maneuver} \cite{SmSt02}. 
\begin{Not} The player in turn to move is \emph{the next player} and the 
other player is \emph{the previous player}.
\end{Not}

\begin{Def}\label{Def:2} Let $k\in \N$. 
In the game of \emph{Blocking-$k$ Wythoff Nim}, denoted by $\W^k$, the 
options are defined as in Wythoff Nim, Definition \ref{Def:1}. 
But before the next player moves, the previous player may declare 
at most $k - 1$ of them as forbidden. When the next player has moved, 
any blocking maneuver is forgotten and has no further impact on the game. 
\end{Def}
Notice that for $k = 1$ this game is Wythoff Nim. 
A player who is unable to move, because all options are forbidden, loses. 

We adapt the standard terminology of $P$- and $N$-positions---the previous 
and next player winning positions respectively---of non-blocking impartial 
games to `$k$-blocking' ditto. 
\begin{Def}\label{Def:3}
The value of (a position of) $W^k$ is $P$ if (strictly) fewer than $k$ of its 
options are $P$, otherwise it is $N$. Denote by $\P _k$ the set 
of $P$-positions of $\W^k$. 
\end{Def}

\begin{figure}[ht]
\centering
\vspace{0.3 cm}
\includegraphics[width=0.815\textwidth]{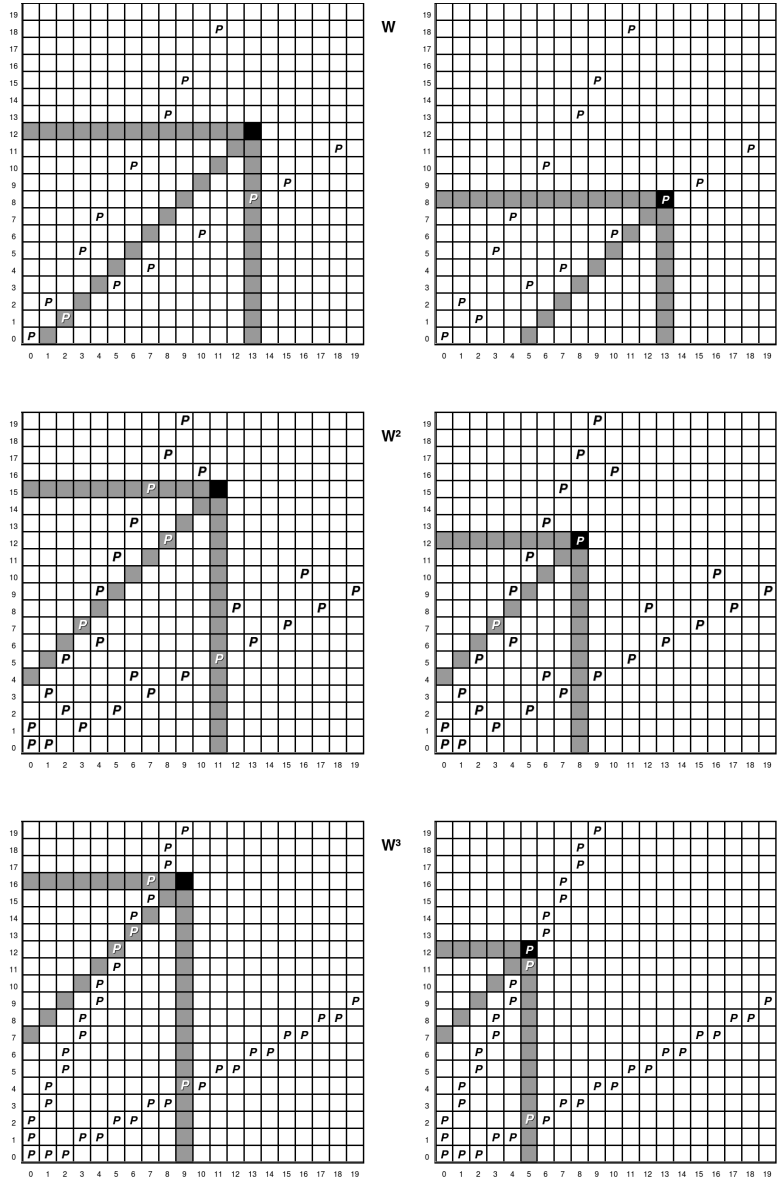}
\caption{
The two figures at the top illustrate options of two instances of Wythoff Nim 
together with its initial $P$-positions.
The middle and lower couples of figures represent 
$\W^2$ and $\W^3$ respectively.
%The figures illustrate two instances of the game 
%$\W^2$ together with the first few $P$-positions. 
For example in the middle left figure the `gray' shaded positions are the 
options of the `black' $N$-position $(11,15)$. 
This position is $N$ since, by rule of game, only one of the two $P$-positions 
in its set of options can be forbidden. 
In contrast, the position $(8,12)$ is $P$ (middle-right) 
since there is precisely one single $P$-position in its set of options. 
It can (and will) be forbidden.}\label{figure:1}
\end{figure}

By this definition, the next player wins if and only if 
the position is $N$. It leads to a recursive definition of the set of 
$P$-positions of $\W^k$, see also Proposition \ref{prop:Wk} on 
page \pageref{prop:Wk}.
%Denote with $\P (G)$ ($\n(G)$) the set of 
%all $P$-positions ($N$-positions) of a game $G$. Let us state our 
%two main results. 
Since both the Wythoff Nim type moves and the blocking maneuvers 
are `symmetric' on the game board it follows 
that the sets of $P$- an $N$-positions are also `symmetric'. 
Hence, the following notation.
\begin{Not}
We use the `symmetric' notation $\{x, y\}$ for unordered pairs of integers, 
that is whenever $(x, y)$ and $(y, x)$ are considered the same.
\end{Not}
Let us explain the main results of this paper, see also Figure \ref{figure:1}.
\begin{Def}\label{Def:R123} Let $\phi = \frac{1 + \sqrt{5}}{2}$ denote 
the Golden ratio. Then 
\begin{align*}
\r_1 &:= \{\{\left \lfloor \phi n \right \rfloor, \left \lfloor \phi^2 n \right \rfloor\}\mid n\in \M\},\\ 
\r_2 &:= \{(0,0)\}\cup \{\{n, 2n + 1\} \mid n\in \M\}\cup
\{(2x + 2, 2y + 2)\mid (x, y)\in \r_{1} \}, \intertext{and}
\r_{3} &:= \{(0,0)\}\cup \{\{n, 2n + 1\}, \{n, 2n + 2\} \mid n\in \M \}\}.
\end{align*}
\end{Def}

\begin{Thm}\label{Thm:1}
The sets $\P_i = \r_i$, $i=1,2,3$. 
\end{Thm}

It is well known that the set $\P_1 = \r_1$ \cite{Wyt07}.
We prove the latter two results in Section \ref{Sec:2}.  Admittedly, 
it surprised me that the solutions of these two problems were so pliable.

In Section \ref{Sec:3} we give a table of a conjectured asymptotic 
`behavior' of $\P_k$ for each $k\in \{4,5,\ldots , 20\}$ and also give a 
brief discussion on a certain family of `Comply games'---in particular 
we define the game $\W_k$ with its set of $P$-positions identical to the set 
of $N$-positions of $\W^k$. 
\subsection{Some general results}
The set $\r_1$ has some frequently studied properties. Namely, the sequences 
$(\left \lfloor \phi n \right \rfloor)$ and 
$(\left \lfloor \phi^2 n \right \rfloor)$ are so-called \emph{complementary 
sequences} of $\N$, e.g. \cite{Fra82}, that is they partition $\N$. 
(This follows from the well known `Beatty's theorem' \cite{Bea}.) In 
this paper we make use of a generalization of this 
concept---often used in the study of so called `(exact) covers by Beatty 
sequences' e.g. \cite{Fra73, Gra73, Heg1}.

\begin{Def}\label{Def:compl}
Let $p\in \N$. Suppose that $A$ is a set of a finite number 
of sequences of non-negative integers. Then $A$ is a \emph{$p$-cover} 
(\emph{cover} if $p = 1$) of another set, say $S\subset \M$, if, 
for each $x\in S$, the total number, $\xi(A,S,x)$, of 
occurrences of $x$, in the sequences of $A$, exceeds or equals $p$. Further, 
$A$ is an \emph{exact} $p$-cover of $S$ if, for all $x$, $\xi(A,S,x) = p$.
\end{Def}

The special case of $S = \N$, $\#A = 2$ and $p = 1$ in this 
definition is `complementarity'. For general $p$ and with $\#A = 2$ the 
term $p$-complementarity is used in \cite{Lar1}. 

Before we move on to the proof of Theorem \ref{Thm:1}, 
let us give some basic results valid for general $\W^k$. 

\begin{Prop}\label{prop:Wk}
Let $k\in \N$ and define $\{\{a_i, b_i\}\mid i\in \M \} = \P_k$, 
where, for all $i$, $a_i\le b_i$ and the ordered pairs $(a_i, b_i)$ 
are in lexicographic order, 
that is $(a_i)$ is non-decreasing and $a_i = a_j$ together with $i < j$ imply 
$b_i < b_j$. Then, 
\begin{enumerate}[(i)]
\item the $0^{th}$ column contains precisely $k$ $P$-positions, namely 
$$(0,0),(0,1),\ldots , (0,k-1),$$ 
\item the set $\{(a_i), (b_i)\mid i\ge k\}$ is an exact $k$-cover of $\N$,
\item for all $d\in \M$, $\#\{i\in \M \mid b_i - a_i = d\}\le k$.
\end{enumerate}
\end{Prop}
%Note that in (i) the first index of the $b$-sequence is 1. The technical 
%reason for this is that we do not want to `double count' 
%the $P$-position $(0,0)$. The first column must obviously contain 
%precisely $k$ $P$-positions, which gives. 
\noindent{\bf Proof.} The case $k = 1$ follows from  well known results on 
Wythoff Nim \cite{Wyt07}. 
%We have already seen 
%that $(0,0)$ is the only $P$-position in the first column, 
%that $(a_i)_{i\in \N}$ and $(b_i)_{i\in \N}$ are complementary and that 
%each `diagonal' of the form $\{(x,x)\mid x\in \M\}$ 
%intersects precisely one $P$-position, which is (iii). 
Hence, let $k > 1$. The item (i) is obvious (see also  (\ref{terminal})). 
For (ii) suppose that there is a least $x'\in \N$ such that 
$$r = \#(\{i\mid a_i = x' \}\cup \{i\mid b_i = x'\})\ne k.$$ Clearly, by 
the blocking rule, this forces $r < k$ for otherwise 
there must trivially exist a non-blocked Nim-type move 
$\boldsymbol x\rightarrow \boldsymbol y$, 
where both $\boldsymbol x, \boldsymbol y\in \P_k$. Suppose that $y$ is the 
largest integer such that $(x', y)\in \P_k$. Then, by the 
blocking rule, for all integers 
\begin{align}\label{zy}
z > y, 
\end{align}
there must exist a $P$-position in the set 
of horizontal and diagonal options of $(x', z)$. (For 
otherwise all $P$-positions 
in the set of options of $(x', z)$ could be blocked off.) But, 
by assumption, the total number of $P$-positions in the columns 
$0, 1, \ldots , x' - 1$ is precisely 
$k(x' - 1)$ and each such position is an option of precisely two positions 
in column $x'$, which contradicts (\ref{zy}). 
Item (iii) is obvious by Definition \ref{Def:2}. 
\hfill $\Box$\\

\begin{Not}
A position (of $\W^k$) is \emph{terminal} if all options may be blocked off by 
the previous player.
\end{Not}

A player who moves to a terminal position may, by Definition \ref{Def:2}, 
be declared the winner. Let $k\in \N$. The terminal positions of $\W^k$ 
are given by the following result. We omit the elementary proof.
%which is easily verified by induction and 
%elementary algebra, so we omit the proof.  
\begin{Prop}
Let $k\in \N$. The set of terminal positions of $\W^k$ is precisely 
\begin{align}\label{terminal}
\T(k) := \{\{x, y\} \mid x\le y < k - 2x, x, y \in \M \}.
\end{align}
The set $\T(k)$ is a \emph{lower ideal}, that is 
$(x, y)\in \T(k)$ implies $(x-i,y-j)\in \T(k)$, 
for all $i\in \{0,1,\ldots ,x\}$ and all $j\in\{0,1,\ldots ,y\}$. 
The number of positions in this set is
\[ \#\T(k) := \left\{ \begin{array}{lll}
3(m + 1)^2 - 2(m + 1) & \mbox{if $k = 3m + 1$,}\\
3(m + 1)^2 & \mbox{if $k = 3m + 2$,}\\
3(m + 1)² + 2(m + 1) & \mbox{if $k = 3(m + 1)$},\end{array} \right. \] 
$m\in \M$.%\hfill $\Box$ 
\end{Prop}

In particular, the set of terminal positions of $\W^2$ and $\W^3$ are 
$\T(2)=\{(0,0),\{0,1\}\}$ ($\#\T(2)=3$) and 
$\T^3=\{(0,0),\{0,1\},\{0,2\}\}$ ($\#\T(3)=5$) respectively.

\section{Proof of the main result}\label{Sec:2}
%We use an `arithmetic-geometric' approach for the proofs in this section. 
%These options may also be grouped in other ways, for example 
%all `Nim-type' options (n) may be taken to denote 
%the horizontal and vertical ones together. 
Given a blocking parameter $k=2$ or $3$ and a position $(x, y)$, we want to 
count the total number of options contained in our candidate 
set of $P$-positions $\r_2$ or $\r_3$ respectively. 
Then we may derive the value of 
$(x, y)$ as follows. The previous player will win if and only if the 
total number of options in the candidate set is strictly less than $k$. With 
this plan in mind, let us define some functions, counting the number of 
options in some specific `candidate set' and of the specific types, 
(v), (d) and (h) respectively.
%(As claimed in Theorem 1.1 and 1.2, this plan has proved succesful 
%for $k=2$ and 3, with candidate sets $\r_2$ and $\r_3$, respectively.)
\begin{Def}\label{Def:count} 
Let $(x, y)\in \B$. Given a set $S\subset \B$, let us define 
\begin{align*}
v_{x, y} = v_{x, y}(S) &:= \#(\{(w,y)\mid x > w \ge 0 \}\cap S),\\
d_{x, y} = d_{x, y}(S) &:= \#(\{(w, z)\mid x - w = y - z > 0 \}\cap S),\\ 
h_{x, y} = h_{x, y}(S) &:= \#(\{(x,z)\mid y > z \ge 0 \}\cap S),\\
f_{x, y} = f_{x, y}(S) &:= d_{x, y}+v_{x, y}+h_{x, y},\\
%n_{x, y} &:= f_{x, y} - d_{x, y} = h_{x, y} + v_{x, y}.
\end{align*}
$w,z\in \M$.
\end{Def}

%Clearly, all functions in Definition \ref{Def:count} are independent of 
%the blocking parameter $k$. In contrast, the value of a 
%position (Definition \ref{Def:3}) obviously depends on $k$. Since all 
%move- and blocking-rules of $\W^k $ are symmetric, the above functions
%  are symmetric. For example, for all $x,y\in \M$, we have that 
%$f_{x, y}=f_{x, y}$ and, for all $x, y\in \M$, the value of $(x, y)$ 
%equals that of $(y, x)$.

\begin{Not}\label{Def:move}
We use the notation $(x_1, x_2) \rightarrow (y_1, y_2)$ if there 
is a Wythoff Nim (Definition \ref{Def:1}) type move from $(x_1, x_2)$ 
to $(y_1, y_2)$. 
%In case there is such a move, it may also be denoted by   
%$(x_1, x_2)\ominus (y_1, y_2) := (x_1 - y_1, x_2 - y_2)$.
\end{Not}
%\begin{Def}
%Let $(x_i)$ and $(y_i)$ be two sequences of integers with the same index set. 
%Then the \emph{difference sequence} of the sequence of pairs $((x_i, y_i))$ 
%is $(y_i - x_i)$.
%\end{Def}
%As we will see in the next section, item (ii) in 
%Proposition \ref{prop:Wk} may, in general, be difficult to improve upon. 
%In particular the difference sequence of 
%$((a_i, b_i))_{i\in \M}$, $(b_i - a_i)_{i\in \M}$, is in general 
%not $k$-complementary. (Of course, the one major exception is 
%for Wythoff Nim, that is, if $k = 1$ this sequence is complementary. 
%See also \cite{HeLa06} for related questions.)

\subsection{Proof of Theorem \ref{Thm:1}.}
With notation as in Definition \ref{Def:R123} and \ref{Def:count}, 
put $S = \r_2$ and let $k = 2$. Hence, we consider the game $\W^2$. 
Then, by the blocking rules in 
Definition \ref{Def:2}, each $P$-position has the property 
that at most one of its options is $P$ and each $N$-position has the 
property that at least two $P$-position are in its set of options.
Thus, the theorem holds if we can prove that the 
value of $(x, y)\in \B$ is $P$ if and only if $f_{x, y}(\r_2) \le 1$. 
Hence notice that (see also (\ref{terminal})) 
\begin{itemize}
\item $f_{0,0} < f_{0, 1} = 1$ and $(0, 0), \{0, 1\}$ are $P$,
\item $x\ge 2$, $f_{0, x} \ge 2$ and $\{0, x\}$ is $N$,
\item $f_{1, 1} = 3$ and $(1, 1)$ is $N$,
\item $f_{1, 2} = 2$ and $\{1, 2\}$ is $N$.
\end{itemize}
Further, the \lq least' $P$-position which is not terminal is $(2, 2)$, namely 
$f_{2, 2} = 1$ since, by the above items, the only option which is 
a $P$-position is $(0, 0)$. 

We divide the rest of the proof of the 
strategy of $\W^2$ into two `classes' depending on 
whether $(x,y)\in \B$ belongs to $ \r_2$ or not.

%\noindent ``$P\rightarrow N$'': 
Suppose that $(x, y)\in \r_2$. That is, we have to prove 
that $f_{x, y}(\r_2)\le 1$. We are done with the cases 
$(x,y) = (0,0), (0,1)$ and $(2,2)$.  We may assume that $1\le x\le y$. \\ 
%Fix a position $(x, y)\in \mathcal{P}$, $x > 0$. \\

\noindent{\it Case 1:} Suppose that $y = 2x + 1$. 
Then, we claim that $h_{x, y} = 0$, $d_{x, y} = 0$ and $v_{x, y} = 1$. \\

\noindent{\it Proof.} The horizontal options 
of $(x, 2x+1)$ are of the form $(z,2x+1)$ with $z<x$. But 
all positions $(r,s)$ in $\r_2$ satisfy 
\begin{align}\label{s2r1}
s\le 2r+1. 
\end{align}
This gives $h_{x, y} = 0$.

The diagonal options are of the form $(z, x+z+1)$, with $0\le z < x$. 
Again, by (\ref{s2r1}), this gives  $d_{x, y} = 0$.

For the vertical options, if $x \le 2$, we are done, so suppose $x > 2$. 
Then, we may use that $\{(2\lfloor \phi n\rfloor + 2)_{n\in \N}, 
(2\lfloor \phi^2 n\rfloor + 2)_{n\in \N}, (2n + 1)_{n\in \N}\}$ is an exact 
cover of $\{3,4,5,\ldots \}$. 

Namely, if $x := 2z + 1$ is odd, we have that 
$$y = 2x + 1 > x = 2z + 1 > z,$$ so that 
$(x, y)\rightarrow (2z+1, z)\in \r_2$. Since $x$ is odd, 
any vertical option in $\r_2$ has to be of this form. 
%is a legal move. 

If, on the other hand, $x := 2z \ge 2$ is even, then, since 
(by \cite{Wyt07}) $(\lfloor \phi n\rfloor )_{n\in \N}$ and 
$(\lfloor \phi^2 n \rfloor )_{n\in \N}$ are complementary, there is 
precisely one $n$ such that either $z = \lfloor \phi n\rfloor + 1$ 
or $ z = \lfloor \phi^2 n \rfloor + 1 $. For the first case 
\begin{align*}
(x, 2x + 1) &= (2z, 4z+1)\\
&= (2 \lfloor\phi n\rfloor + 2, 4 \lfloor \phi^2 n\rfloor + 3)\rightarrow (2 \lfloor \phi n\rfloor +2,  2\lfloor \phi ^2 n\rfloor + 2)\in \r_2.
\end{align*}
%is a legal move. 
The second case is similar. But, since $x$ is even, any option in $\r_2$ has 
to be precisely of one of these forms. We may conclude that $v_{x,y}=1$.\\

\noindent{\it Case 2:} Suppose that $x = 2\lfloor \phi n \rfloor +2 $ and 
$y = 2\lfloor \phi ^2 n\rfloor + 2$, for some $n\in \M$. 
Then, we claim that $d_{x, y} = 1$ and $v_{x, y} = h_{x, y} = 0$.\\

\noindent {\it Proof.} 
%The difference sequence of $(0,0)$ together with 
%$((x, 2x + 1))$ is $D := (x)_{x\in \M} $. 
%In particular, each even difference occurs precisely once in $D$, as in 
%$d := (2\lfloor \phi^2 n \rfloor - 2\lfloor \phi n \rfloor)_{n\in \N} = (2n)$. 
If $n = 0$, we are done, so suppose that $n > 0$. We have that 
$$(2\lfloor \phi n \rfloor + 2 , 2\lfloor \phi n + n \rfloor + 2) 
-(2n -1, 4n - 1) 
= (2\lfloor \phi n \rfloor - 2n + 3 , 2\lfloor \phi n \rfloor - 2n + 3)$$ 
is a diagonal move in Wythoff Nim (and where the `-'sign denotes vector 
subtraction). This gives $d_{x,y} \ge 1$. We may partition the differences of 
the coordinates of the positions in $\r_2$ into two sequences, 
$$((2n+1)-n)_{n\in \M} = (n)_{n\in \N}$$ and 
$$(2\lfloor \phi^2 n \rfloor + 2 - (2\lfloor \phi n \rfloor + 2))_{n\in \M} 
= (2n)_{n\in \N}$$ respectively. These sequences are strictly increasing, 
which gives $d_{x,y} = 1$.
 
For the second part we may 
apply the same argument as in Case 1, but in the other direction. 
Namely, 
%$x = 2\lfloor \phi^t n \rfloor + 2$, $t\in \{1, 2\}$, gives 
$2x + 1 \ge 4\lfloor \phi n \rfloor + 3 > 2\lfloor \phi^2 n\rfloor + 2 
> 2\lfloor \phi n\rfloor + 2$, which implies that all nim-type options 
belong to the set $\B\setminus \r_2$. 

We are done with the first class.
%\noindent ``$N\rightarrow P$'': 
Hence assume that $(x, y)\not \in \r_2$. 
 That is, we have to prove that $f_{x, y}(\r_2)\ge~2 $.\\

\noindent{\it Case 3:} Suppose $y > 2x + 1$. 
Then we claim that $v_{x, y} = 2$, $h_{x, y} = 0$ and $d_{x, y} = 0$.\\ 

\noindent {\it Proof.} By the first argument in Case 1, the latter two claims 
are obvious. Notice that the set of sequences 
\{$(n)_{n\in \N}$, $(2n + 1)_{n\in \M}$, 
$(2\lfloor \phi n\rfloor + 2)_{n\in \M }$,  
$(2\lfloor \phi^2 n \rfloor+2 )_{n\in \M}$\} constitute an 
exact 2-cover of $\N$. This gives $v_{x, y} = 2$.\\

\noindent{\it Case 4:} Suppose $0 < x \le y < 2x + 1$. 
Then, we claim that either 
\begin{enumerate}[(i)]
\item $d_{x, y} = 1$ and $h_{x, y} + v_{x,y}\ge 1$, or 
\item $d_{x, y} = 2$.
\end{enumerate}

\noindent {\it Proof.} We consider three cases.
\begin{enumerate}[(a)]
\item $y > \phi x$,
\item $y < \phi x$ and $y - x$ even,
\item $y < \phi x$ and $y - x$ odd.
\end{enumerate}

In case (a), $v_{x, y} = 1$ is verified as in Case 1. For $d_{x, y} = 1$, 
%. This means that there is a 
%main-diagonal Wythoff-type move $(x, y)\rightarrow (n, 2n+1)$ for some $n < x$.
it suffices to demonstrate that $(x, y)-(z, 2z + 1) = (x-z, y-2z-1)$, 
is a legal diagonal move for some $z\in \M$. 
Thus, it suffices to prove that $x - z = y - 2z - 1$ holds 
together with $0 < z < x$ and $2z + 1 < y$. But this follows since the 
definition of $y$ implies $z+1 = y-x \le 2x-x= x$.

In case (b) we get $d_{x,y} = 2$ by 
$(2\lfloor \phi^2 n\rfloor  - 2\lfloor \phi n \rfloor )_{n\in \N} = (2n)$ 
and an analog reasoning as in the latter part of (a). (Hence this is (ii).)

In case (c) we may again use the latter 
argument in (a), but, for parity reasons, 
there are no diagonal options of the first type in (b), 
so we need to return to case (i) and thus verify that 
$h_{x, y} + v_{x,y}\ge 1$. Since $y-x$ is odd we get that 
precisely one of $x$ or $y$ 
must be of the form $2z + 1$, $z\in \M$.
Suppose that $z < x = 2z+1 < y$. Then $(x, y)\rightarrow (2z + 1, z)$ 
gives $v_{x, y}\ge 1$. If, on the other hand, $x\le  y = 2z + 1 < \phi x$,  
then $(x, y)\rightarrow (z, 2z + 1)$ is legal since 
$z < \frac{\phi x - 1}{2}<x$, which gives $h_{x,y}\ge 1$.
%For the first part of item (i) it suffices to use the 
%same argument as in the proof of Case 2. For the second part, the first 
%part of the proof of Case 1 together with $y >  \phi x $ 
%gives $v_{x, y}\ge 1$. If either $x$ or $y$ is odd, a 
%similar argument gives either $v_{x, y}\ge 1$ or  $v_{x, y}\ge 1$, 
%so assume both $x$ and $y$ are even with $x\le y <  \phi x $. 
%Then $y - x = 2n$ for some $n\in \N$. We claim that 
%$(x, x + 2n)\rightarrow (2\lceil \phi n\rceil , 2\lceil \phi n + n\rceil )$. 
%We have that $x + 2n <  \phi x $. This 
%gives $2n < (\phi - 1) x  $. Multiplying both sides 
%with $\phi$ gives $2n\phi < x$. Since $x = 2t$ is even we may simplify 
%to $\lceil n\phi \rceil \le t$ since $t$ is an integer. But the case 
%$x = 2\lceil n\phi \rceil$ is already excluded since $(x, x + 2n)\in \mathcal{N}$. This gives the claim and the theorem. 

We are done with $\W^2$'s part of the proof. 
Therefore, let $S = \r_3$, $k=3$ and consider the game $\W^3$. 
Then one needs to prove that $(x,y)\in \r_3$, $x\le y$, if and only 
if $f_{x, y}(\r_3)\le 2$. Suppose that $(x, y)\in \r_3$ with $x\le y$. 
Then we claim that $d_{x, y}\le 1$, $h_{x, y} = 0$ and $v_{x, y}\le d_{x, y} + 1$. 
Otherwise, if $(x, y)\not \in \r_3$ and $y > 2x + 2$, then we claim 
that $v_{x, y} = 3$, or, if $y < 2x + 1$, then we claim that 
$h_{x, y} = v_{x, y} = 1$ and  $d_{x, y} \ge 1$. Each case is almost immediate by 
definition of $\r_3$ and Figure \ref{figure:1}, so we omit further details. 
\hfill $\Box$

\section{Discussion}\label{Sec:3}
%In \cite[Remark 4,5]{Lar3} the concept of an \emph{$l$-fold split} 
%is introduced. 
One obvious direction of future research is to try and 
classify the $P$-positions of the games $\W^k$, $k\ge 4$. Let $k\in \N$ 
and, as in Proposition \ref{prop:Wk}, let $\{\{a_i, b_i\}\mid i\in \M \}=\P_k$ 
denote the set of $P$-positions of $\W^k$ (with, for all $i$, $b_i\ge a_i$). 
Then we ask if 
\begin{align}\label{nosplit}
\lim_{i\rightarrow \infty}\frac{b_{i}}{a_{i}} 
\end{align}
exists. If not, then we wonder if the set of aggregation points of 
$\lim_{i\rightarrow \infty}\frac{b_{i}}{a_{i}}$ is finite. More precisely: 
For some (least) $2 \le l = l(k)\in \N$, does 
there exist $l$ sequences $t^j$, $j\in\{1, 2, \ldots , l\}$, such that 
pairwise distinct asymptotic limits 
\begin{align}\label{split}
\lim_{i\rightarrow \infty}\frac{b_{t^j_i}}{a_{t^j_i}} 
\end{align}
exist? We conjecture that, for each $k\in \{4,5, \ldots , 20\}$ 
in Table \ref{t1}, $l(k)$ is given by the number of entries in row $k$. 
%The cases $k = 1, 2$ 
%and $3$ are discussed in Section \ref{Sec:1} and \ref{Sec:2}.

In \cite{Lar2} another generalization of Wythoff Nim is studied, namely the 
family of Generalized Diagonal Wythoff Nim games and a 
so-called \emph{split} of sequences 
of ordered pairs is defined. In particular a sequence of pairs 
$((a_i, b_i))$ is said to \emph{split} 
if (\ref{nosplit}) is not satisfied but (\ref{split}) 
is (for some $l\ge 2$). In that paper one conjectures quite 
remarkable asymptotic `splits' of $P$-positions for certain 
games---supported by numerous computer 
simulations and figures. However, the only proof of a `splitting' of 
$P$-positions given in that paper is 
the much weaker statement that (\ref{nosplit}) is not satisfied, and 
it is only given for one 
particular game called $(1,2)$GDWN---a game which extends the diagonal options 
of Wythoff Nim and also allows moves of the types 
$(x,y)\rightarrow (x-i,y-2i)$, $x\ge i>0, y\ge 2i>0$ and 
$(x,y)\rightarrow(x-2i,y-i)$, $x\ge 2i>0, y\ge i>0$. 
\begin{Rem}
In this paper we have proved that the `upper' (above the main diagonal)
$P$-positions of $\W^2$ split. I am not aware of any other such result, 
of a split of `the upper' $P$-positions of an impartial game, 
in particular not on a variation of Wythoff Nim. (if we drop the `upper' 
condition then one may obviously regard the $P$-positions of, for example, 
Wythoff Nim as a splitting sequence, see \cite{Lar2}).
\end{Rem}
At the end of this section, we provide tables of the 
first few $P$-positions for $\W^4$, $\W^5$ and $\W^6$ respectively. 
As an appetizer for future research on Blocking Wythoff Nim games, 
let us motivate the conjectured asymptote 
of (\ref{nosplit}) in Table \ref{t1} for the case $\W^4$, that is that  
$\lim_{i\rightarrow \infty}\frac{b_{i}}{a_{i}}$ 
exists and equals $\sqrt{2}+1$. If (\ref{nosplit}) holds, with 
$\lim a_i/i = \alpha$ and $\lim b_i/i = \beta$ real numbers, then, by 
Proposition \ref{prop:Wk}, also $\alpha^{-1} + \beta^{-1} = 4$ holds. 
Also, by Table \ref{t4}, one hypothesis is that 
$\delta_n = b_n-a_n = n/2 + O(1)$, where $O(1)$ denotes some 
bounded function. This  gives $\beta - \alpha = 1/2$ and so, by
elementary algebra we get 
$$\lim_{i\rightarrow \infty}\frac{b_{i}}{a_{i}}= \frac{\beta}{\alpha} =\sqrt{2} +1.$$

\subsection{Comply- versus blocking-games}
Let $k\in \N$ and let $G_k$ denote the following 
`comply'-variation of any impartial game $G$. The 
previous player is requested to propose at least $k$ of the options 
of $G$ as allowed next-player options in $G_k$ (and these are all options). 
After the next player has moved, this `comply-maneuver' is forgotten and 
has no further impact on the game. 
The player who is unable to propose at least $k$ next-player options is the 
loser\footnote{Our `comply' games provide a subtle variation 
to those in \cite{SmSt02}.}. 
It is not hard to check that this game has the `reverse' strategy 
of $G^k$, where $G^k$ is defined in analogy with $\W^k$ 
(in Definition \ref{Def:2} exchange Wythoff Nim for $G$). By this we mean that 
the set of $P$-positions of $G_k$ is precisely the set 
of $N$-positions of $G^k$. 

%The first non-trivial example one might come 
%to think of, to apply this comply-maneuver, is the game of one pile Nim 
%(with the trivial blocking maneuver, that is $k=1$), 
 
Thus, as an example, let Nim$_1$ denote 
the comply-variation of Nim where the previous player has to propose 
at least one Nim option. In this game the empty pile 
is $N$ (the previous player loses because he cannot propose 
any option). Each non-empty pile is $P$, 
since the previous player will propose the empty pile as the only 
available option for the next player. 
Recall that the only $P$-position of Nim (without blocking maneuver) is 
the empty pile.

This discussion motivates why we, in the definition of $\W^k$ 
(Definition \ref{Def:2}), let the previous player forbid $k - 1$, rather 
than $k$ options. To propose at least $k$ options is the 
`complement' of forbidding fewer than $k$ options---and it is not 
a big surprise that the set of $P$-positions of $\W^k$ are `complementary' 
to those of $\W_k$. (Another more `algorithmic' way of thinking of this 
choice of notation is that (the position of) $\W^k$ a priory belongs 
to the set of forbidden options.)

Other blocking maneuvers on Wythoff Nim have been studied in the literature. 
Let $k\in \N$. 
In the game of $k$-blocking Wythoff Nim \cite{HeLa06, Lar09, FrPe} the 
blocking maneuver only constrains moves of type (d) in 
Definition \ref{Def:1}. Otherwise the rules are the same as in this paper. 
For the purpose of this section, denote this game by $\W^k$N.
In another variation, the game of 
Wythoff $k$-blocking Nim \cite{Lar1}, the blocking maneuver only constrains 
moves of type (h) or (v). Denote this game by $\W\text{N}^k$. 
(Both these game families are actually defined, and solved, as restrictions 
of $m$-Wythoff Nim, \cite{Fra82}).

The above discussion leads us to the following round up of this paper. 
Namely, let us attempt to define the corresponding `comply rules' 
of $\W^k$N and $\W\text{N}^k$, that is we want to find rules of games $\W_k$N 
and $\W\text{N}_k$ such that the $P$-positions of $\W\text{N}^k$ 
correspond precisely 
to the $N$-positions of $\W\text{N}_k$ and the $P$-positions of 
$\W^k$N correspond precisely to the $N$-positions of $\W_k$N. We claim 
that the correct comply maneuver of $\W_k$N is: Propose at least $k$ options 
of the type (d) \emph{or} at least one 
Nim-type, (h) or (v), option. Similarly, for the game $\W\text{N}_k$ 
the comply rules are: Propose at least $k$ Nim-type options \emph{or} 
at least one option of type (d). Otherwise the rules are as in $\W_k$. 

Why are these the correct rules? Suppose, for example, 
that $(x, y)$ is $P$ in $\W^k$N, say. 
Then there are at most $k-1$ (d) type $P$-positions and no 
Nim-type $P$-position at all in the set of options of $(x,y)$. We have to 
demonstrate that $(x,y)$ is $N$ in $\W_k$N, that is that the previous 
player cannot propose $k$ (d)-type positions, all of them $N$, neither 
can he propose a single Nim-type $N$-position. But, by symmetry, this claim 
follows from a straightforward inductive argument, assuming 
complementarity of the $P$-positions in the two games on the $<x$ 
indexed columns. We omit further details. I could think of a few ways 
to continue this discussion, but it is nicer to leave some of 
this territory wide open in favor of others ideas.

\begin{table}[!ht]
\begin{center}
\begin{tabular}{ l || c | c | c | c | }
%\hline
%\backslashbox{k}{l} &1 & 2 & 3 & 4& \hline
  $k=$ & $l=1$ & $l=2$  & $l=3$ & $l=4$ \\ \hline
  1   &$\phi$ & &&\\ \hline
  2   &$\phi$ &2 && \\\hline  
  3   &2 & && \\\hline
  4   &$1+\sqrt{2}$ & && \\\hline
  5   &1.476 &2.5 && \\\hline
  6   &1.28 &2.0 &2.5& \\\hline
  7   & 2.5 & && \\\hline
  8   & 2.0 & 2.5 && \\\hline
  9   & 1.34&2.5 && \\\hline
  10  &  1.59&2.0 &2.5& \\\hline
  11  &2.5 & && \\\hline
  12  & 1.74&2.5 && \\\hline  
  13  & 1.2&2.5 && \\\hline 
  14  &  2.5& && \\\hline
  15  &  $1+\sqrt{2}$&2.58 && \\\hline
  16  & 1.426 &2.5 &2.6& \\\hline
  17  & 1.12 &2.0 &2.5&2.6 \\\hline  
  18  &2.35 &2.5 &2.6& \\\hline
  19  &1.88 &2.5 &2.6&\\\hline
  20  &1.28 &2.5 &2.6&\\\hline

\end{tabular}
\end{center}\caption{The entries in this table are the estimated/conjectured 
quotients $\lim_{i\rightarrow \infty}\frac{b_{t^j_i}}{a_{t^j_i}}$ for 
$j\in\{1,2,\ldots ,l\}$ and the respective game $\W^k$. The cases $k = 2, 3$ 
are resolved in Theorem \ref{Thm:1} and $k = 1$ is Wythoff Nim, 
$\phi=\frac{\sqrt{5}+1}{2}$.}\label{t1}
\end{table}

\begin{table}[ht!]
\begin{center}
\begin{tabular}{ l | c | c | c|| l | c | c | c || l | c | c | c }
$n$ & $a_n$ & $b_n$ & $\delta_n$ & $n$ & $a_n$ & $b_n$ & $\delta_n$ & $n$ & $a_n$ & $b_n$ & $\delta_n$ \\ \hline 
0 & 0 & 0 & 0 &30 & 10 & 25 & 15 & 60 & 20 & 50 & 30 \\ 
1 & 0 & 1 & 1 &31 & 10 & 26 & 16 & 61 & 21 & 51 & 30 \\
2 & 0 & 2 & 2 & 32 & 11 & 26 & 15 & 62 & 21 & 52 & 31 \\
3 & 0 & 3 & 3 & 33 & 11 & 27 & 16 & 63 & 22 & 53 & 31 \\
4 & 1 & 1 & 0 & 34 & 11 & 28 & 17 & 64 & 22 & 54 & 32 \\ 
5 & 1 & 4 & 3 & 35 & 12 & 29 & 17 & 65 & 22 & 55 & 33 \\
6 & 1 & 5 & 4 & 36 & 12 & 30 & 18 & 66 & 23 & 55 & 32 \\ 
7 & 2 & 3 & 1 & 37 & 12 & 31 & 19 & 67 & 23 & 56 & 33 \\ 
8 & 2 & 6 & 4 & 38 & 13 & 31 & 18 & 68 & 23 & 57 & 34 \\ 
9 & 2 & 7 & 5 & 39 & 13 & 32 & 19 & 69 & 24 & 58 & 34 \\ 
10 & 3 & 8 & 5 & 40 & 13 & 33 & 20 & 70 & 24 & 59 & 35 \\ 
11 & 3 & 9 & 6 & 41 & 14 & 34 & 20 & 71 & 24 & 60 & 36 \\ 
12 & 4 & 6 & 2 & 42 & 14 & 35 & 21 & 72 & 25 & 60 & 35 \\ 
13 & 4 & 10 & 6 & 43 & 14 & 36 & 22 & 73 & 25 & 61 & 36 \\ 
14 & 4 & 11 & 7 & 44 & 15 & 37 & 22 & 74 & 25 & 62 & 37 \\ 
15 & 5 & 12 & 7 & 45 & 15 & 38 & 23 & 75 & 26 & 63 & 37 \\ 
16 & 5 & 13 & 8 & 46 & 16 & 37 & 21 & 76 & 26 & 64 & 38 \\ 
17 & 5 & 14 & 9 & 47 & 16 & 39 & 23 & 77 & 27 & 65 & 38 \\ 
18 & 6 & 15 & 9 & 48 & 16 & 40 & 24 & 78 & 27 & 66 & 39 \\ 
19 & 6 & 16 & 10 & 49 & 17 & 41 & 24 & 79 & 27 & 67 & 40 \\ 
20 & 7 & 15 & 8 & 50 & 17 & 42 & 25 & 80 & 28 & 67 & 39 \\ 
21 & 7 & 17 & 10 & 51 & 17 & 43 & 26 & 81 & 28 & 68 & 40 \\ 
22 & 7 & 18 & 11 & 52 & 18 & 43 & 25 & 82 & 28 & 69 & 41 \\ 
23 & 8 & 19 & 11 & 53 & 18 & 44 & 26 & 83 & 29 & 70 & 41 \\ 
24 & 8 & 20 & 12 & 54 & 18 & 45 & 27 & 84 & 29 & 71 & 42 \\ 
25 & 8 & 21 & 13 & 55 & 19 & 46 & 27 & 85 & 29 & 72 & 43 \\ 
26 & 9 & 21 & 12 & 56 & 19 & 47 & 28 & 86 & 30 & 72 & 42 \\ 
27 & 9 & 22 & 13 & 57 & 19 & 48 & 29 & 87 & 30 & 73 & 43 \\ 
28 & 9 & 23 & 14 & 58 & 20 & 48 & 28 & 88 & 30 & 74 & 44 \\ 
29 & 10 & 24 & 14 & 59 & 20 & 49 & 29 & 89 & 31 & 75 & 44 \\ 
\end{tabular}
\end{center}\caption{The first few $P$-positions of $\W^4$, $\{a_n, b_n\}$, and 
the corresponding differences $\delta_n:=b_n - a_n$.}\label{t4}
\end{table}
\clearpage
\begin{table}[ht!]
\begin{center}
\begin{tabular}{ l | c | c | c|| l | c | c | c || l | c | c | c }
$n$ & $a_n$ & $b_n$ & $\delta_n$ & $n$ & $a_n$ & $b_n$ & $\delta_n$ & $n$ & $a_n$ & $b_n$ & $\delta_n$ \\ \hline 
0 & 0 & 0 & 0 & 30 & 8 & 21 & 13 & 60 & 17 & 24 & 7 \\
1 & 0 & 1 & 1 & 31 & 8 & 22 & 14 & 61 & 17 & 43 & 26 \\ 
2 & 0 & 2 & 2 & 32 & 8 & 23 & 15 & 62 & 17 & 44 & 27 \\ 
3 & 0 & 3 & 3 & 33 & 9 & 23 & 14 & 63 & 17 & 45 & 28 \\ 
4 & 0 & 4 & 4 & 34 & 9 & 24 & 15 & 64 & 18 & 46 & 28 \\ 
5 & 1 & 1 & 0 & 35 & 9 & 25 & 16 & 65 & 18 & 47 & 29 \\ 
6 & 1 & 2 & 1 & 36 & 10 & 14 & 4 & 66 & 18 & 48 & 30 \\ 
7 & 1 & 5 & 4 & 37 & 10 & 26 & 16 & 67 & 19 & 25 & 6 \\ 
8 & 1 & 6 & 5 & 38 & 10 & 27 & 17 & 68 & 19 & 48 & 29 \\ 
9 & 2 & 4 & 2 & 39 & 10 & 28 & 18 & 69 & 19 & 49 & 30 \\ 
10 & 2 & 7 & 5 & 40 & 11 & 16 & 5 & 70 & 19 & 50 & 31 \\ 
11 & 2 & 8 & 6 & 41 & 11 & 28 & 17 & 71 & 20 & 29 & 9 \\ 
12 & 3 & 3 & 0 & 42 & 11 & 29 & 18 & 72 & 20 & 51 & 31 \\ 
13 & 3 & 6 & 3 & 43 & 11 & 30 & 19 & 73 & 20 & 52 & 32 \\ 
14 & 3 & 9 & 6 & 44 & 12 & 15 & 3 & 74 & 20 & 53 & 33 \\ 
15 & 3 & 10 & 7 & 45 & 12 & 31 & 19 & 75 & 21 & 31 & 10 \\  
16 & 4 & 11 & 7 & 46 & 12 & 32 & 20 & 76 & 21 & 53 & 32 \\ 
17 & 4 & 12 & 8 & 47 & 12 & 33 & 21 & 77 & 21 & 54 & 33 \\  
18 & 4 & 13 & 9 & 48 & 13 & 33 & 20 & 78 & 21 & 55 & 34 \\ 
19 & 5 & 7 & 2 & 49 & 13 & 34 & 21 & 79 & 22 & 30 & 8 \\ 
20 & 5 & 13 & 8 & 50 & 13 & 35 & 22 & 80 & 22 & 56 & 34 \\  
21 & 5 & 14 & 9 & 51 & 14 & 36 & 22 & 81 & 22 & 57 & 35 \\ 
22 & 5 & 15 & 10 & 52 & 14 & 37 & 23 & 82 & 22 & 58 & 36 \\ 
23 & 6 & 16 & 10 & 53 & 14 & 38 & 24 & 83 & 23 & 58 & 35 \\  
24 & 6 & 17 & 11 & 54 & 15 & 38 & 23 & 84 & 23 & 59 & 36 \\ 
25 & 6 & 18 & 12 & 55 & 15 & 39 & 24 & 85 & 23 & 60 & 37 \\
26 & 7 & 18 & 11 & 56 & 15 & 40 & 25 & 86 & 24 & 61 & 37 \\ 
27 & 7 & 19 & 12 & 57 & 16 & 41 & 25 & 87 & 24 & 62 & 38 \\ 
28 & 7 & 20 & 13 & 58 & 16 & 42 & 26 & 88 & 24 & 63 & 39 \\ 
29 & 8 & 9 & 1 & 59 & 16 & 43 & 27 & 89 & 25 & 63 & 38 \\ 
\end{tabular}
\end{center}\caption{The first few $P$-positions of $\W^5$, $\{a_n, b_n\}$, and 
the corresponding differences $\delta_n := b_n - a_n$.}
\end{table}
\clearpage
\begin{table}[ht!]
\begin{center}
\begin{tabular}{ l | c | c | c|| l | c | c | c || l | c | c | c }
$n$ & $a_n$ & $b_n$ & $\delta_n$ & $n$ & $a_n$ & $b_n$ & $\delta_n$ & $n$ & $a_n$ & $b_n$ & $\delta_n$ \\ \hline 
0 & 0 & 0 & 0 & 30 & 6 & 19 & 13 & 60 & 13 & 36 & 23 \\ 
1 & 0 & 1 & 1 & 31 & 7 & 15 & 8 & 61 & 14 & 27 & 13 \\ 
2 & 0 & 2 & 2 & 32 & 7 & 19 & 12 & 62 & 14 & 37 & 23 \\ 
3 & 0 & 3 & 3 & 33 & 7 & 20 & 13 & 63 & 14 & 38 & 24 \\ 
4 & 0 & 4 & 4 & 34 & 7 & 21 & 14 & 64 & 14 & 39 & 25 \\ 
5 & 0 & 5 & 5 & 35 & 8 & 9 & 1 & 65 & 15 & 30 & 15 \\ 
6 & 1 & 1 & 0 & 36 & 8 & 22 & 14 & 66 & 15 & 39 & 24 \\ 
7 & 1 & 2 & 1 & 37 & 8 & 23 & 15 & 67 & 15 & 40 & 25 \\ 
8 & 1 & 3 & 2 & 38 & 8 & 24 & 16 & 68 & 15 & 41 & 26  \\
9 & 1 & 6 & 5 & 39 & 9 & 16 & 7 & 69 & 16 & 32 & 16  \\
10 & 1 & 7 & 6 & 40 & 9 & 24 & 15 & 70 & 16 & 42 & 26  \\
11 & 2 & 4 & 2 & 41 & 9 & 25 & 16 & 71 & 16 & 43 & 27  \\
12 & 2 & 5 & 3 & 42 & 9 & 26 & 17 & 72 & 16 & 44 & 28   \\
13 & 2 & 8 & 6 & 43 & 10 & 20 & 10 & 73 & 17 & 18 & 1   \\
14 & 2 & 9 & 7 & 44 & 10 & 27 & 17 & 74 & 17 & 31 & 14  \\ 
15 & 3 & 6 & 3 & 45 & 10 & 28 & 18 & 75 & 17 & 44 & 27  \\ 
16 & 3 & 7 & 4 & 46 & 10 & 29 & 19 & 76 & 17 & 45 & 28  \\ 
17 & 3 & 10 & 7 & 47 & 11 & 11 & 0 & 77 & 17 & 46 & 29  \\ 
18 & 3 & 11 & 8 & 48 & 11 & 22 & 11 & 78 & 18 & 35 & 17  \\
19 & 4 & 8 & 4 & 49 & 11 & 29 & 18 & 79 & 18 & 47 & 29   \\
20 & 4 & 12 & 8 & 50 & 11 & 30 & 19 & 80 & 18 & 48 & 30  \\ 
21 & 4 & 13 & 9 & 51 & 11 & 31 & 20 & 81 & 18 & 49 & 31   \\
22 & 4 & 14 & 10 & 52 & 12 & 21 & 9 & 82 & 19 & 37 & 18   \\
23 & 5 & 10 & 5 & 53 & 12 & 32 & 20 & 83 & 19 & 49 & 30   \\
24 & 5 & 14 & 9 & 54 & 12 & 33 & 21 & 84 & 19 & 50 & 31   \\
25 & 5 & 15 & 10 & 55 & 12 & 34 & 22 & 85 & 19 & 51 & 32   \\
26 & 5 & 16 & 11 & 56 & 13 & 13 & 0 & 86 & 20 & 40 & 20  \\
27 & 6 & 12 & 6 & 57 & 13 & 25 & 12 & 87 & 20 & 52 & 32  \\
28 & 6 & 17 & 11 & 58 & 13 & 34 & 21 & 88 & 20 & 53 & 33  \\
29 & 6 & 18 & 12 & 59 & 13 & 35 & 22 & 89 & 20 & 54 & 34  \\
\end{tabular}
\end{center}\caption{The first few $P$-positions of $\W^6$, $\{a_n, b_n\}$, and 
the corresponding differences $\delta_n := b_n - a_n$.}
\end{table}
%\section{From $p$-complementarity to complementarity}
%Let $A = \{(a_i), (b_i), (c_i) \text{ and } (d_i)\mid i\in \N \}$ 
%denote a set of four sequences of positive integers with 
%the following properties: 
%\begin{itemize}
%\item $A$ is a $k$-cover of $\N$,
%\item $\alpha := \lim_{i\in \N}\frac{b_i}{a_i}$ exists,
%\item $\beta := \lim_{i\in \N}\frac{d_i}{c_i}$ exists,
%\item $\alpha - \beta = \delta > 0$.
%\end{itemize}
%Let us denote by \emph{the $k$-split problem}, the question 
%to find an impartial game  whose set of $P$-positions is 
%$\{\{a_i, b_i\},\{c_i, d_i\}\mid i\in \N\}\cup \{\boldsymbol 0\}.$

\end{document}